# Optimal non-adaptive approximation of convex bodies by polytopes

Kamenev G.K.[1][0000-0002-1809-2017]

[1] Federal Research Center of Informatics and Management of the RAS, Vavilova str, 40, Moscow, Russian Federation.
gkk@ccas.ru

**Abstract.** In this paper we consider the problem of constructing numerical algorithms for approximating of convex compact bodies in $d$-dimensional Euclidean space by polytopes with any given accuracy. It is well known that optimal with respect to the order algorithms produces polytopes for which the accuracy in Hausdorff metric is inversely proportional to the number of vertices (faces) in the degree of 2/($d$-1). Numerical approximation algorithms can be adaptive (active) when the vertices or faces are constructed successively, depending on the information obtained in the process of approximation, and non-adaptive (passive) when parameters of algorithms are defined on the basis of a priory information available. Approximation algorithms differ in the use of operations applied to the approximated body. Most common are indicator, support and distance (Minkowski) functions calculations. Some optimal active algorithms for arbitrary bodies approximation are known using support or distance function calculation operation. Optimal passive algorithms for smooth bodies approximation are known using support function calculation operation and extremal curvature information. It is known that there are no optimal non-adaptive algorithms for arbitrary bodies approximation using support function calculation operation. We consider optimal non-adaptive algorithms for arbitrary bodies approximation using projection function calculation operation.

## 1    The problem of polyhedral non-adaptive approximation

The problem of approximation of convex compact bodies by polyhedra is a classical one. Many fundamental results are known [1, 2]. The development of real-life algorithms is of practical importance especially in operational research problems for example in multicriteria decision-making [3, 4]. The best possible approximation of convex compact bodies by polytopes with a given number of vertices or hyperfaces is well studied theoretically [1, 2]. Adaptive (active) and non-adaptive (passive) algorithms of polyhedral approximation differ. Both use a priori information about the approximated body, but in adaptive algorithms, it is specified by the information obtained in the process of approximation and taken into account in it. For non-adaptive (passive) algorithms, parameters are defined only on a priory information available.

Most modern adaptive and non-adaptive algorithms of polyhedral approximation assume the possibility of calculating the support function (the maximum of the linear form) on the approximated body. There are a number of optimal adaptive algorithms [1, 2, 5, 6]. Unfortunately, they are non-stable for multidimensional case in the presence of calculation errors and special tools of topological combinatorics are used (see [3], 6.1) for their real-life implementation making them very complicated in coding. In the presence of a priori information about the smoothness of the approximated body and the maximum radius of curvature of the surface both adaptive and non-adaptive approaches gives the same optimal order of the convergence rate on the number of hyperfaces of the approximating polyhedron [7, 8]. In the non-smooth case, as well as in the absence of information on the smoothness or on the maximum radius of curvature of the approximated body, the order of convergence rate by the number of hyperfaces of non-adaptive algorithms based on the calculation of the support function is significantly lower [9]. In the classical results concerning optimal polyhedral approximation of convex bodies [10, 11] the technique for approximating of arbitrary convex bodies based on the operation of projection of external point can be found. We implements this technique to construct simple in realization non-adaptive algorithms based on the operation of projection (finding the minimum of a quadratic form), which have an optimal order of convergence rate by the number of hyperfaces.

## 2 Optimal non-adaptive methods

### 2.1 Optimal approximation methods

Consider the Euclidean space $E^d$ with the scalar product $<\cdot,\cdot>$, distance $\rho(\cdot,\cdot)$ and the norm $\|\cdot\|$. Let $B_r$ denote a ball of radius $r$ and $B$ be a unit ball with center at the origin, $S$ – be a sphere (hypersphere) of directions, i.e. $\partial B$. Let $C \subset E^d$ – be a convex compact body and $P \subset E^d$ – be a compact convex polytope. Using $m^t(P)$ we denote the number of vertices, and using $m^f(P)$ we denote the number of hyperfaces (faces of maximum dimension). For definiteness, we assume that $C \subset B$. For convex bodies $C_1$ and $C_2$, consider the Hausdorff metric

$$\delta(C_1, C_2) = \max\{\sup\{\rho(x, C_2): x \in C_1\}, \sup\{\rho(x, C_1): x \in C_2\}\}.$$

It is proved [1, 2] that, for some constant $c$, there exist such polytopes $P_m$ such that for any $m = m^t(P_m)$ (or $m = m^f(P_m)$) it holds true that

$$\delta(C, P_m) \leq c / m^{2/(d-1)},$$

and for bodies with twice continuously differentiable boundary these estimates can not be improved. The method of polyhedral approximation, which allows to construct approximating polytopes $P$ with the property

$$\delta(C, P) \leq c / [m^f(P)]^{2/(d-1)}, \; m^f(P) \to \infty,$$

will be called the *optimal method* (more precisely, optimal in order) *by the number of hyperfaces*.

## 2.2 Non-adaptive approximation methods based on the calculation of the support function

For $u \in S$ we introduce the notation of the support function, the support hyperplane, and the support half-space:

$$g(u, C) = \max \{<u, x>: x \in C\},$$

$$l(u, C) = \{x \in E^d: <u, x> = g(u, C)\},$$

$$L(u, C) = \{x \in E^d: <u, x> \leq g(u, C)\}.$$

For an arbitrary $p \in E^d$ we will use the notations $g(u, p) = <u, p>$, $l(u, p) = \{x \in E^d: <u, x> = <u, p>\}$, $L(u, p) = \{x \in E^d: <u, x> \leq <u, p>\}$, which does not contradict the previously introduced definitions of the support function, hyperplane and half-space.

For a given $U \subset S$ we denote

$$P(C, U) = \cap \{L(u, C): u \in U\} = \cap \{ \{x \in E^d: <u, x> \leq g(u, C)\} : u \in U \}.$$

Let $U(m)$ be a subset of $S$ consisting of $m$ elements. It is obvious that if $P(C, U)$ is a polytope then

$$m^f(P(C, U(m))) \leq m.$$

Through $\delta(U(m))$ we denote *the approximation accuracy guaranteed by the set $U(m)$*:

$$\delta(U(m)) = \max\{\delta(C, P(C, U(m))): C \subset B\}.$$

We will say that a *non-adaptive approximation method based on the calculation of the support function* is defined if for any $m$ the set $U(m)$ is given such that $\delta(U(m)) \to 0$, $m \to \infty$. There are non-adaptive approximation methods based on the calculation of the support function of convex compact bodies. For example in [6] it is proved, that for any $m$ it holds true that

$$\delta(U(m)) \leq c / m^{1/(d-1)}$$

at some constant $c$. However, as shown in [9], for any method of selection of $U(m)$ there is a constant $c^*$ such that for each $m$ there is $C \subset B$, for which

$$\delta(C, P(C, U(m))) \geq c^* / [m^f(P(C, U(m)))]^{1/(d-1)},$$

following

$$\delta(U(m)) \geq c^* / m^{1/(d-1)}.$$

Thus, *non-adaptive methods of convex bodies approximation based on the calculation of its supporting function are not optimal by the number of hyperfaces*.

## 3  Projection technique

The set *U* is called the $\varepsilon$-net of some set if any point of the set is located from *U* at a distance not greater than $\varepsilon$. For $p \notin C$ we denote by proj($p$, $C$) the projection of point $p$ on *C*. the point proj($p$, $C$) is the only one and is the solution of the quadratic minimization problem:

$$\text{proj}(p, C) = \arg\min\{\rho(p, x): x \in C\}, p \notin C.$$

We denote $F: E^d \setminus C \to S$,

$$F(p) = (p - \text{proj}(p, C)) / \|p - \text{proj}(p, C)\|.$$

In [11] it is shown that using the $\varepsilon$-nets of the ball $B_2$, consisting of *m* points, and projecting to the surface of the body $C \subset B$ it is possible to construct a polyhedron with the property

$$\delta(C, P) \leq c / [m^f(P)]^{2/(d-1)}.$$

In [10] a method for constructing a polyhedron with the property

$$\delta(C, P) \leq c / [m^t(P)]^{2/(d-1)}.$$

We will use the basic ideas from [11] to develop our own technique for constructing an optimal non-adaptive method based on the design operation.
Let $\beta > 0$ and $U \subset \partial B_{(1+\beta)}$ are given. We denote

$$P^+(C, U) = P(C, F(U)) = \cap\{L(u, C): u \in F(U)\}.$$

To construct $P^+(C, U)$ it is necessary to find $p^+ = \text{proj}(p, C)$ for each $p \in U$ (i.e. to solve the problem of bilinear programming) and put

$$u = (p - p^+) / \|p - p^+\|, g(u, C) = <u, p^+>.$$

**Theorem 1.** *Let $\beta > 0$ is given and U is $\varepsilon$-net $\partial B_{(1+\beta)}$, $\varepsilon < \beta$. Then*

$$\delta(C, P^+(C, U)) \leq \varepsilon^2 / (\beta^2 - \varepsilon^2)^{1/2}.$$

Note that the factor $1/(\beta^2 - \varepsilon^2)^{1/2}$ can be chosen arbitrarily close to $1/\beta$ when choosing a sufficiently small $\varepsilon$ compared to $\beta$.
*Proof of Theorem 1.* Let $C \subset B$ and $P = P^+(C, U)$. Let

$$p \in \partial P: \rho(p, C) = \delta(C, P^+(C, U)).$$

Let $x = \text{proj}(p, C)$ and $p^{**} \in \partial B_{(1+\beta)}$, such that $x = \text{proj}(p^{**}, C)$. Since *U* is $\varepsilon$-net of $\partial B_{(1+\beta)}$, then there exists $q^{**} \in U$, such that $\rho(p^{**}, q^{**}) \leq \varepsilon$. We denote

$$y = \text{proj}(q^{**}, C), u = F(p^{**}) = (p^{**} - x) / \|p^{**} - x\|, v = F(q^{**}) = (q^{**} - y) / \|q^{**} - x\|.$$

Let $p^*=x+\beta u$, $q^*=y+\beta v$. Since $x$ and $y$ belong to $\partial C$, then $p^*$ и $q^*$ belong to $\partial(C+B_\beta)$. Projecting from the outside on a convex set does not increase the distance, so $\rho(p^*, q^*) \leq \varepsilon$ and $\rho(x, y) \leq \varepsilon$. We have

$$\|(x+\beta u) - (y+\beta v)\| = \rho(p^*, q^*) \leq \varepsilon.$$

Then

$$\max\{\|x-y\|, \beta\|u-v\|\} \leq \varepsilon.$$

Let $\alpha$ be the angle between $u$ and $v$. Since $\sin \alpha \leq \|u - v\|$, then $\sin \alpha \leq \varepsilon/\beta$. Let $z$ be the projection of $x$ on $l = l(u, x) \cap l(v, y)$ and $m$ be the projection of $x$ on $l(v, y)$. The projection of $m$ on $l$ is $z$ (otherwise, $z$ would not be the projection of $x$ since

$$\rho(z, x) = (\rho(z, m)^2 + \rho(x, m)^2)^{1/2},$$

and the value of angle $\angle xzm$ is $\alpha$. Let $p^{***}$ be the intersection of the ray $[x, p^{**})$ with $l(v, y)$. It is obvious that $p \in [x, p^{***}]$. Therefore

$$\delta(C, P^+(C, U)) = \rho(x, p) \leq \rho(x, p^{***}) = \rho(x, m) / \cos \alpha =$$

$$= \rho(x, z) \sin \alpha / \cos \alpha \leq \rho(x, y) \sin \alpha / \cos \alpha \leq \varepsilon^2 / \beta [1-(\varepsilon/\beta)^2)]^{1/2},$$

following the statement of the theorem..

## 4  Constructing efficient $\varepsilon$-nets of the hypersphere

### 4.1  Efficient covering of the hypersphere

Let us now consider the possibility of constructing $\varepsilon$-nets on the surface of the sphere. We denote by $m(\varepsilon)$ the minimum number of points of the $\varepsilon$-net of $S$. For a given $\varepsilon$ we know quite little [12]. However, the asymptotic value of

$$m(\varepsilon) \approx \vartheta/\varepsilon^{d-1}, \; \varepsilon \to 0,$$

where $\vartheta = \vartheta_{d-1}(d\pi_d/\pi_{d-1})$. Here $\vartheta_d$ denotes the minimum density of the coverage of space $E^d$ by balls of fixed radius [12], and $\pi_d = \pi^{d/2}/\Gamma((d/2)+1)$ denotes the volume of a single ball, $\Gamma$ is gamma function. Precisely known only values $\vartheta_1=1$, $\vartheta_2=2\pi\sqrt(27)$. The construction method of $\varepsilon$-nets $U(m)$, $m \to \infty$, such that there is a constant $\Theta$, for which

$$m \leq \Theta/\varepsilon^{d-1},$$

we will call the method of building an *efficient covering of the hypersphere* (with constant $\Theta$). The value $\eta = (\vartheta/\Theta)^{(d-1)}$ we call the asymptotic efficiency of this method. For the same number of points of the $\varepsilon$-net, this value shows how many times worse, asymptotically, this method allows us to cover the surface of the hypersphere, in comparison with the optimal one.

### 4.2 Methods based on lattice coverings of Euclidian space

The most common in applications is the use of a uniform grid (lattice) in polar coordinates. In this case, it is true [13]

$$M \leq \Theta_1 / \varepsilon^{d-1}, m \to \infty, \Theta_1 = 2[\pi^2(d-1)/4]^{(d-1)/2}.$$

The asymptotic efficiency of this method for large $d$ does not exceed $2e/(\pi\sqrt{d})$ [13]. Thus, the efficiency of the method of $\varepsilon$-net constructing using a uniform grid in polar coordinates decreases rapidly with the dimension of the problem. Another disadvantage of this approach is the need for a sharp jump to increase the number of points of the $\varepsilon$-net ($m$ from $k^{d-1}$ to $(k+1)^{d-1}$, where $k$ is the number of divisions coordinate) if necessary to increase the accuracy of the coverage, while changing all points. This same problem applies to other methods of $\varepsilon$-nets construction based on lattice coverings of Euclidian space.

### 4.3 Methods based on the approximation of the balls by polytopes

To construct the efficient coverings of the hypersphere we can use efficient methods of polyhedral approximation of the ball [14, 15]. They may be based on the following relationship.

**Theorem 2.** *Let $Q$ be the polytope inscribed in the ball $B$, then the set of its vertices will be $\varepsilon$-net of $S$, where*

$$\varepsilon = (2\delta(Q, B))^{1/2}.$$

The proof of this theorem is completely analogous to the proof of Property 4 [16] or theorem 2 [17] for the internal hypersphere metric.

There are efficient adaptive methods for polyhedral approximation of the convex bodies [6]. Among them, the most studied theoretically and in numerical experiments is "Estimate Refinement Method" [1, 15]. In this method, the sequence of points monotonically increasing in the number $m$ is constructed on the surface of the sphere. In this sequence each initial part forms a suboptimal $\varepsilon(m)$-net. Though the algorithm of the "Estimate Refinement Method" is rather difficult to implement, such a sequence is calculated once and can be further applied directly in non-adaptive methods. Appropriate $\varepsilon(m)$-nets are given in [18] for dimensions 3 to 10. As shown in [14, 15], the asymptotic efficiency of the sequence of nets is not below

$$(\vartheta_{d-1}/\delta_{d-1})^{1/(d-1)}/2 \geq [2(d-1)/d]^{1/2}/2,$$

where $\delta_d$ denotes the maximum packing density of the space $E^d$ by balls of a fixed radius [12]. Thus, the lower estimate of the asymptotic efficiency of this method with increasing dimension approaches $\sqrt{2}/2$.

# 5 Optimal non-adaptive approximation methods based on the calculation of the projection function

## 5.1 A class of non-adaptive methods of polyhedral approximation based on the operation of projecting and effective covering of the hypersphere

Now let a certain method of hypersphere covering is given. We will say that a *non-adaptive approximation method based on the calculation of the projection function* is defined if for any $m$ the set $U(m)$ is given such that $\delta(C, P^+(C, U(m))) \to 0$, $m \to \infty$. We will construct the class of optimal non-adaptive methods of polyhedral approximation based on the operation of projecting and effective covering of the hypersphere.

**Lemma 1.** Let $\beta > 0$ is given and $U(m)$ be $\varepsilon(m)$-net on $S$. Then if $(1+\beta)\varepsilon(m) < \beta$ it holds

$$\delta(C, P^+(C, (1+\beta)U(m))) \leq (1+\beta)^2 \varepsilon(m)^2/(\beta^2 - [(1+\beta)\varepsilon(m)]^2)^{1/2}.$$

*Proof of the Lemma 1.* Since $U(m)$ is $\varepsilon(m)$-net on $S$, then $(1+\beta)U(m)$ is $(1+\beta)\varepsilon(m)$-net on $\partial B_{(1+\beta)}$. So from theorem 1 it follows the assertion of lemma under condition $(1+\beta)\varepsilon(m) < \beta$

**Lemma 2.** Let $\beta > 0$ is given and $U(m)$ be $\varepsilon(m)$-net on $S$ with property $m \leq \Theta/\varepsilon^{d-1}$, for some constant $\Theta$. Then for any $\lambda > 0$ there exists $m(\lambda)$, such that for $m \geq m(\lambda)$ it holds

$$\delta(C, P^+(C, (1+\beta)U(m))) \leq (1+\lambda)A(\beta)/m^{2/(d-1)}, A(\beta) = (1+\beta)^2 \Theta^{2/(d-1)}/\beta.$$

*Proof of the Lemma 2.* The value $\varepsilon(m)$ can be chosen to be arbitrarily small, from which we obtain the statement of the Lemma 2 from Lemma 1.

Now let us suppose that $U(m)$ is given such that $U(m)$, $m \to \infty$, is efficient covering of the hypersphere.

**Theorem 3.** *Let $P^m = P^+(C, 2U(m))$, $m \to \infty$, be a sequence of polytopes generated by the non-adaptive polyhedral approximation method based on the projection operation and on constructing efficient coverings $U(m)$ of the sphere with the constant $\Theta$. Then for any $\lambda > 0$ there exists $m(\lambda)$, such that for $m \geq m(\lambda)$ it holds*

$$\delta(C, P^m) \leq (1+\lambda)A/m^{2/(d-1)}, A = 4\Theta^{2/(d-1)}.$$

*Proof of theorem 3.* It is sufficient to note that the minimum value of $A(\beta)$ in Lemma 2 is achieved when $\beta = 1$ and is equal to $A = 4\Theta^{2/(d-1)}$.

## 5.2 Optimal non-adaptive methods of polyhedral approximation based on the operation of projecting and effective covering of the hypersphere

In conclusion, we give a detailed description of the method.

Let us construct an approximation of a convex compact body $C \subset B$ with accuracy $\delta$ in the Hausdorff metric. Then it is necessary:

— to choose $\varepsilon < 1/2$: $4\varepsilon^2/(1-4\varepsilon^2)^{1/2} \leq \delta$;

- to build $\varepsilon$-net $U$ of $S$ using the existing method of effective covering of the hypersphere;
- to find for each point $u \in U$ the point

$$v = F(2u) = (2u - \text{proj}(2u, C)) / \|2u - \text{proj}(2u, C)\|$$

by solving the quadratic minimization problem on $C$, and to define

$$g(v, C) = <v, \text{proj}(2u, C)>.$$

- to take as a solution a polyhedron

$$P^+(C, 2U) = \cap \{ \{x \in E^d: <v, x> \leq g(v, C)\}: v \in F(2U) \}.$$

**Corollary.** Non-adaptive method for polyhedral approximation based on the projection operation and efficient covering of the hypersphere is optimal on order of the number of hypefaces of the approximating polytope.

The statement of the Corollary follows from theorem 3 taking in account that $m^f(P^+(C, U(m))) \leq m$.

This method is easily generalized to the case of convex compact bodies without the restriction $C \subset B$. To do this, it is enough to find an external ball for the body and to scale the parameters of the method to its radius.

This work was supported by the Russian Foundation for Basic Research, project no. 18-01-00465a.